\documentclass[aos,preprint]{imsart}

\RequirePackage{amsthm,amsmath,amsfonts,amssymb,mathrsfs}
\RequirePackage[round,authoryear]{natbib}
\RequirePackage[colorlinks,citecolor=blue,urlcolor=blue]{hyperref}
\RequirePackage{graphicx,subfigure}
\usepackage{bbm,bm}
\usepackage{enumerate}
\usepackage{caption}
\usepackage{xcolor}
\usepackage{algorithm}
\usepackage{algpseudocode}
\usepackage{booktabs} 
\usepackage{graphicx}

\startlocaldefs
\theoremstyle{plain}

\newtheorem{thm}{Theorem}[section]
\newtheorem{lemma}[thm]{Lemma}
\newtheorem{remark}{Remark}[section]

\newtheorem{cor}{Corollary}[section]
\newtheorem{proposition}{Proposition}[section]
\theoremstyle{remark}
\newtheorem{definition}[thm]{Definition}
\newtheorem{example}{Example}[section]

\endlocaldefs

\newcommand{\bbU}{{\bf U}}

\newcommand{\bbA}{{\bf A}}

\newcommand{\bbX}{{\bf X}}

\newcommand{\bbZ}{{\bf Z}}

\newcommand{\bSi}{\pmb \Sigma}

\newcommand{\bgL}{{\bf \Lambda}}

\newcommand{\bPi}{{\bf \Pi}}

\newcommand{\bbI}{{\bf I}}

\newcommand{\bqn}{\begin{eqnarray*}}
	\newcommand{\eqn}{\end{eqnarray*}}

\newcommand{\bqa}{\begin{eqnarray}}
	\newcommand{\eqa}{\end{eqnarray}}

\newcommand{\al}{\alpha}

\newcommand{\CYRS}{\CYRS}

\begin{document}
	
	\begin{frontmatter}
		\title{High-dimensional ridgeless least squares interpolation under spiked covariance structures}
		\runtitle{}
		
		\begin{aug}
			\author[1]{\fnms{Zhijun} \snm{Liu}\ead[label=e1]
				{liuzj@mail.neu.edu.cn}} and
			\author[2]{\fnms{Dandan} \snm{Jiang}\thanks{Corresponding author.}\ead[label=e2]{jiangdd@xjtu.edu.cn}} 
	    	
			\address[1]{College of Sciences, Northeastern University, China.
			\printead{e1}
			}
			\address[2]{School of Mathematics and Statistics, Xi’an Jiao Tong University, China.
			\printead{e2}
			}
		\end{aug}
		
\begin{abstract}
This paper investigates the asymptotic behavior of the out-of-sample prediction risk of the high-dimensional ridgeless least-squares estimator when the feature dimension $p$ and the sample size $n$ grow proportionally. We consider a generalized spiked population covariance model with multiple latent factors, where the number of spiked eigenvalues may remain finite or increase with $n$, and the spiked eigenvalues may be bounded or diverge at arbitrary rates. 
Beyond characterizing the impact of covariance spectra, we reveal a new mechanism underlying benign overfitting: the prediction behavior of ridgeless interpolation is fundamentally governed by the alignment between the regression coefficient  $\bm\beta$ and the spiked eigenspaces of the population covariance matrix. In particular, we show that the signal energy distributed along latent spike directions determines whether interpolation leads to benign, tempered, or catastrophic overfitting.
Our theoretical framework establishes sharp prediction risk limits under minimal moment conditions, requiring only finite fourth moments rather than Gaussianity. We characterize how the number, strength, and geometric structure of the spikes jointly influence the double-descent phenomenon. These results provide a unified understanding of when latent covariance structures facilitate or hinder generalization in overparameterized regression.

\end{abstract}
		
		
	\end{frontmatter}

\section{Introduction}
The modern practice of machine learning has challenged the classical statistical view of model complexity. According to the traditional bias--variance trade-off, increasing model complexity first reduces bias but eventually increases variance, leading to a U-shaped test error curve. However, many contemporary overparametrized models, including neural networks and high-dimensional linear predictors, are capable of interpolating the training data while still achieving strong out-of-sample performance. This apparent contradiction has motivated a new line of theoretical work on the generalization behavior of interpolating estimators.

A central concept in this direction is the double descent phenomenon. \cite{Belkin19} proposed a unified risk curve that extends the classical U-shaped bias--variance trade-off beyond the interpolation threshold. In this picture, the prediction risk may increase as the model approaches the point where it first interpolates the training data, often exhibiting a peak near the interpolation threshold, but then decreases again as the model becomes increasingly overparametrized. This framework provides a conceptual explanation for why highly expressive models can generalize well despite having enough capacity to fit the training data exactly.

Subsequent studies have developed more precise mathematical models for this phenomenon. In particular, \cite{Belkin20} studied two models of double descent for weak features and analyzed the risk of least-squares and least-norm predictors. Their results demonstrate that, when additional features are weak but informative, the test risk can peak around the interpolation boundary and then decrease in the overparametrized regime. These studies indicate that double descent is not merely an empirical artifact of complex learning systems, but can arise even in relatively simple linear regression models. Building on this perspective, \cite{Hastie22} provided a systematic asymptotic analysis of ridgeless least squares and ridge regression in proportional-dimensional linear models, thereby offering a precise statistical characterization of interpolation and generalization in high dimensions.

Most existing theoretical analyses, however, are conducted under relatively restrictive assumptions regarding the population covariance structure. In particular, \cite{Hastie22} studied the isotropic setting, in which the population covariance matrix is the identity matrix, and also considered model misspecification. Although these settings reveal many important features of high-dimensional regression, they do not fully capture the more general non-isotropic covariance structures commonly encountered in applications. In factor models, signal-plus-noise models, and related high-dimensional statistical settings, the population covariance matrix often exhibits several prominent eigenvalues associated with latent factors or signal directions, while the remaining eigenvalues constitute a bulk noise component. These considerations motivate the study of prediction risk under a generalized latent spiked covariance model.

In this paper, we study the asymptotic prediction risk of the min-norm least squares estimator under a generalized spiked population covariance structure. Specifically, we assume independent and identically distributed (i.i.d.) data $(y_i,\bm x_i)$, $i\in \{1,\dots,n\}$, where $\bm x_i\in \mathbb{R}^p$ is a feature vector and $y_i\in \mathbb{R}$ is a response variable.
The observations are generated according to the linear model
\begin{align}
	\label{model1}	\left(\bm x_i, \epsilon_i\right) & \sim P_x \times P_\epsilon, \quad i=1, \ldots, n, \\
	\label{model2}	y_i & =\bm x_i^{\top} \bm{\beta}+\epsilon_i, \quad i=1, \ldots, n,
\end{align}
where $P_x$ is a distribution on $\mathbb{R}^p$ such that $\mathbb{E}\left(\bm x_i\right)=0, \operatorname{Cov}\left(\bm x_i\right)=\bSi$, and $\bSi$ is a spiked model (see details in Assumption \ref{ass1}), and $P_\epsilon$ is a distribution on $\mathbb{R}$ such that $\mathbb{E}\left(\epsilon_i\right)=0, \operatorname{Var}\left(\epsilon_i\right)=\sigma^2$. 
Our population covariance matrix $\bSi$ follows a generalized spiked model,
\begin{align}\label{spikejiegou}
	\bSi=\bbU\operatorname{diag}(\alpha_1,\ldots,\alpha_M,1,\ldots,1)\bbU^\top,
\end{align}
where $\bbU$ is an orthogonal matrix and the spiked eigenvalues $\{\alpha_1,\ldots,\alpha_M\}$ represent a few strong population directions.
We work in the proportional asymptotic regime $(p/n\to\gamma\in(0,\infty)),$ allowing both the underparameterized case $(\gamma<1)$ and the overparameterized case $(\gamma>1).$ The main focus is the overparameterized setting, where the least squares objective has infinitely many solutions and the min-norm least squares estimator, also known as the ridgeless least squares estimator, selects the solution with minimum Euclidean norm. Specifically, it can be formulated as:
we collect the responses in a vector $\bm y=(y_1,y_2,\dots,y_n)\in\mathbb{R}^n$, and the features in a matrix $\bbX=(\bm x_1,\bm x_2,\dots,\bm x_n)^{\top}\in \mathbb{R}^{n\times p}$.
Then the min-norm least squares estimator of $\bm y$ on $\bbX$, is defined by
\begin{align}\label{mininorm}
	\hat{\bm \beta}=\arg \min \left\{\|\bm b\|_2: \bm b \text { minimizes }\|\bm y-\bbX \bm b\|_2^2\right\}.
\end{align}
Alternatively, it can be expressed as $\hat{\bm\beta}=\left(\bbX^\top \bbX\right)^{+} \bbX^\top \bm y$, where $\left(\bbX^\top \bbX\right)^{+}$ is the pseudoinverse of $\bbX^\top \bbX$.

Our population covariance matrix follows a generalized spiked model,
\begin{align}\label{spikejiegou}
\bSi=\bbU\operatorname{diag}(\alpha_1,\ldots,\alpha_M,1,\ldots,1)\bbU^\top,
\end{align}
where $\bbU$ is an orthogonal matrix and the spiked eigenvalues $\{\alpha_1,\ldots,\alpha_M\}$ represent a few strong population directions. The setting (\ref{spikejiegou}) is motivated by the classical spiked model in which a few large eigenvalues of the population covariance matrix are assumed to be well separated from the remaining eigenvalues. \cite{Johnstone01} 
first introduced the spiked population model as the non-null case where all eigenvalues of $\bSi$ are equal to one except for a fixed small number of spikes, i.e.,
$$
\operatorname{Spec}\left(\bSi\right)=\{\alpha_1, \cdots, \alpha_M, \underbrace{1, \cdots, 1}_{p-M}\}.
$$ 
Following Johnstone's development, many efforts have been put into quantifying the effect caused by spiked eigenvalues $\left\{\alpha_k, 1 \leq k \leq M\right\}$ on $M$ extreme sample eigenvalues $\left\{s_k, 1 \leq k \leq M\right\}$. To name a few, under Johnstone's spiked model settings, \cite{baik06} thoroughly studied the almost sure limits of the extreme sample eigenvalues under the Marčenko-Pastur regime when $p, n \rightarrow \infty, p / n \rightarrow \gamma \in(0, \infty)$. They showed that these limits are different when the corresponding population spiked eigenvalues are larger or smaller than critical values $1+\sqrt{\gamma}$ and $1-\sqrt{\gamma}$. 
The spiked model has served as the foundation for a rich theory of principal component analysis through the performance of extreme eigenvalues, as discussed in
\cite{Paul07,Bai08,bai12}. This model extends the classical spiked covariance framework introduced in high-dimensional principal component analysis and provides a natural mathematical formulation for structured population heterogeneity. In a variety of domains, including economics and wireless communications, the spiked model serves as the foundational framework for numerous practical problems; see \cite{bai02,Johnstone17}.

The quantity of central interest is the out-of-sample prediction risk. For a new test point $\bm x_0$, independent of the training data, the conditional prediction risk of an estimator $\hat{\bm\beta}$ is
$$
R_X(\hat{\bm\beta};\bm\beta)=
\mathbb E[
(\bm x_0^\top\hat{\bm\beta}-\bm x_0^\top\bm\beta)^2
\mid \bbX
].
$$
For the min-norm least squares estimator, this risk admits a bias--variance decomposition involving a projection onto the null space of the sample covariance matrix and the trace functional $\operatorname{tr}(\hat{\bSi}^{+}\bSi)$ (see Lemma \ref{lem1} in Section \ref{notation}). 
For isotropic designs with $\bSi=\bbI_p$, existing analyses characterize prediction risk through the spectral properties of sample covariance matrices. In contrast, under a generalized latent spiked covariance model, the prediction risk is governed not only by the limiting spectral distribution of the sample covariance matrix but also by the geometric interaction between the regression signal and the latent covariance structure. In particular, the asymptotic behavior of the empirical eigenvectors and their alignment with the population spiked eigenspaces become essential in determining the bias component of the risk. This coupling determines the bias behavior and plays a fundamental role in distinguishing benign, tempered, and catastrophic overfitting regimes. It introduces a fundamentally new challenge beyond eigenvalue-based analyses.

Several recent works have studied related questions under spiked or signal-plus-noise structures. For instance, \cite{M19} considered a regime in which $p/n\to\infty$ and the population covariance matrix has diverging spiked eigenvalues, and derives bounds for the risk of the min-norm least squares estimator. More recently, \cite{li24} studied a signal-plus-noise model with a single-spike population structure and derived the asymptotic generalization error in the proportional regime under rotational invariance assumptions on the bulk component. \cite{li25} further investigated the effects of spike strength and target-spike alignment, again primarily in a single-spike signal-plus-noise setting with strong invariance assumptions. Compared with these results, our work allows a more general multi-spike population covariance structure and relaxes the distributional requirements on the noise component. In particular, our analysis requires only suitable moment conditions rather than rotational bi-invariance of the bulk noise.

Our main contribution are summarized as follows.
\begin{itemize}
\item First, we derive deterministic asymptotic limits for the prediction risk of ridgeless least squares under a generalized spiked covariance model. Notably, our framework allows the spiked eigenvalues to diverge without imposing any restrictions on their rates of divergence. Moreover, the number of spikes may either remain finite or grow with $n$.
\item Second, we provide a detailed characterization of how the aspect ratio, the number and strengths of the spikes, the signal direction, and the noise level jointly determine generalization performance. In particular, we analyze the role of target--spike alignment, namely, the degree of alignment between $\bm\beta$ and the population spike eigenvectors, in determining the generalization error. These results provide a more comprehensive understanding of double descent and interpolation in linear regression with a structured population covariance matrix. We further conduct extensive numerical simulations to illustrate the theoretical findings and provide additional intuition.
\item Finally, we investigate three overfitting regimes under the spike structure--namely, benign, tempered, and catastrophic overfitting--as well as the behavior of the double descent curve under the generalized spiked model.
\end{itemize}

The technical analysis combines several tools from large-dimensional random matrix theory. The main technical difficulty comes from the interaction between the sample covariance matrix and the population spike directions. In the isotropic case $\bSi=\bbI_p$, the risk can be analyzed mainly through the limiting spectral distribution of the sample covariance matrix. Under a generalized spiked covariance model, however, the leading sample eigenvalues may separate from the bulk spectrum, and the associated sample eigenvectors may have non-negligible asymptotic overlap with the population spike eigenvectors. Consequently, risk functionals such as $\operatorname{tr}(\hat{\bSi}^{+}\bSi)$ and $\bm\beta^\top\bPi\bSi\bPi\bm\beta$ depend not only on eigenvalue limits, but also on eigenvector alignment. To handle this issue, we decompose the relevant matrix functionals into bulk and spike components, and then apply asymptotic results for both eigenvalues and eigenvectors in spiked covariance models. This allows us to identify the contribution of each spike to the limiting prediction risk and to describe how these contributions depend on the alignment between the target vector $\bm\beta$ and the spiked eigenspaces.


\textcolor{black}{The rest of the paper is organized as follows. Section \ref{notation} introduces model assumptions, and technical conditions. Section \ref{mainresult} presents the main asymptotic results for the prediction risk under the generalized spiked covariance model. In Section \ref{stat imp}, we  investigate three overfitting regimes under the spike structure—namely, benign, tempered, and catastrophic overfitting—as well as the behavior of the double descent curve in this setting. 
Some useful lemmas are given in the Appendix.}


Throughout the paper, bold uppercase letters and bold italic lowercase letters are used to denote matrices and vectors, respectively, while scalars are denoted by regular letters. For a matrix $\bbA$, we use
$\operatorname{tr}(\bbA)$, $\bbA^\top$, and $\bbA^+$
to denote its trace, transpose, and Moore--Penrose pseudoinverse, respectively. We use $\bm e_i$ to denote the $i$th standard basis vector, whose components are all zero except that its $i$th component is equal to one. The notation $f'$ denotes the derivative of a function $f$. We write $\lambda_i^{\bbA}$ for the $i$th largest eigenvalue of the matrix $\bbA$.
Throughout the paper, $o(1)$ denotes a deterministic quantity that converges to zero, whereas $o_p(1)$ denotes a random quantity that converges to zero in probability. We write $A_n=\Omega(B_n)$ to mean that $B_n=O(A_n)$. $ \left\langle \bm a,\bm b\right\rangle $ represents the inner product of two vectors $\bm a$ and $\bm b$. Finally, $C$ denotes a generic positive constant whose value may vary from line to line.

\section{Model formulation}\label{notation}
For the min-norm least squares estimator $\bm{\hat{\beta}}$ under the spiked covariance structure for $\bbX$, we first define the out-of-sample prediction risk, and then provide the bias--variance decomposition of the risk.
Consider a test point $\bm x_0 \sim P_x$, independent of the training data. For an estimator $\bm{\hat{\beta}}$ (a function of the training data $\bbX, \bm y$), its out-of-sample prediction risk (or simply the risk) is defined as 
$$R_X(\bm{\hat{\beta}} ; \bm{\beta})=\mathbb{E}[(\bm x_0^{\top} \bm{\hat{\beta}}-\bm x_0^{\top} \bm{\beta})^2 \mid \bbX]=\mathbb{E}[\|\bm{\hat{\beta}}-\bm{\beta}\|_{\bSi}^2 \mid \bbX],$$
where $\|\bm x\|_{\bSi}^2=\bm x^{\top} \bSi \bm x$.  The risk admits the following bias-variance decomposition:
$$
R_X(\bm{\hat{\beta}} ; \bm{\beta})=\underbrace{\|\mathbb{E}(\bm{\hat{\beta}} \mid \bbX)-\bm{\beta}\|_{\bSi}^2}_{B_X(\bm{\hat{\beta}} ; \bm{\beta})}+\underbrace{\operatorname{tr}[\operatorname{Cov}(\bm{\hat{\beta}} \mid \bbX) \bSi]}_{V_X(\bm{\hat{\beta}} ; \bm{\beta})} .
$$
Here $B_X(\bm{\hat{\beta}} ; \bm{\beta})$ is the bias term and $V_X(\bm{\hat{\beta}} ; \bm{\beta})$ is the variance term.  When observations are from model (\ref{model1})
 and (\ref{model2}), equivalent representations of the bias term and variance term are given by the following lemma.
\begin{lemma}[Lemma 1 in \cite{Hastie22}]\label{lem1}
	Under the model (\ref{model1})
	and (\ref{model2}), the min-norm least squares estimator (\ref{mininorm}) has bias and variance
	$$
	B_X(\bm{\hat{\beta}} ; \bm{\beta})=\bm{\beta}^{\top} \bPi \bSi \bPi \bm{\beta} \quad \text { and } \quad V_X(\bm{\hat{\beta}} ; \bm{\beta})=\frac{\sigma^2}{n} \operatorname{tr}(\hat{\bSi}^{+} \bSi),
	$$
	where $\hat{\bSi}=\bbX^{\top} \bbX / n$ is the (uncentered) sample covariance of $\bbX$, and $\bPi=\bbI-\hat{\bSi}^{+} \hat{\bSi}$ is the projection onto the null space of $\bbX$.
\end{lemma}
Throughout the rest of the paper, $\hat{\bSi}$ is defined as in Lemma \ref{lem1}. Before we derive the asymptotic properties of the prediction risk by random matrix theory, we first introduce the fundamental theories in the following. 
\begin{definition}[Empirical spectral distribution]
	For any matrix $\bbA$ with real eigenvalues, the
	empirical spectral distribution of $\bbA$ is denoted by 
	\begin{align*}
		F^{\bbA}\left(x \right)=\frac{1}{p}\left(\text{number of eigenvalues of }  \bbA \leq x \right).  
    \end{align*}
\end{definition}	

\begin{definition}[Stieltjes transform]
	For any function of bounded variation $F$ on the real line, its Stieltjes transform is defined by
	$$
	m_{F}(z)=\int \frac{1}{x-z} \mathrm{~d} F(x), \quad z \in \mathbb{C}^{+} :=\{z \in \mathbb{C}: \Im z>0\}.
	$$
\end{definition}

The assumptions used to obtain the results in this paper are as follows.
	\newtheorem{assumption}{Assumption}[]
\begin{assumption} \label{ass2}
	{ As $\min\{p,n\}\to\infty$, the ratio of the dimension-to-sample size} $ \gamma_{n}:={p}/{n}\rightarrow \gamma>0. $
\end{assumption} 
\begin{assumption} \label{ass1}
For observations $\bm{x_i}$, we assume $\bm{x_i}=\bSi^{1/2}\bm{z_i}$, where	$\bm{z_i}$ are i.i.d. random vectors, and its elements $ \{z_{ij},  1\leq i\leq p,  1\leq j\leq n \} $ have zero means, unit variances, and $\mathbbm{E}|z_{ij}|^4<\infty$. 	For the population covariance matrix, we assume $\bSi=\bbU\bgL\bbU^{\top}$, where $\bbU=(\bm u_1,\bm u_2,\dots,\bm u_p)$ is an orthogonal matrix, and
\begin{align*}
	\bgL=diag(\al_1,\al_2,\dots,\al_M,1,1,\dots,1).
\end{align*}
Here $M$ spikes $\al_1\geq\al_2\geq\dots\geq\al_M>1+\sqrt{\gamma}$. They can be bounded or tend to infinity. $M$ is fixed or $M\rightarrow\infty$ and $M/n^{1/4}\rightarrow0$. 
\end{assumption}
\begin{remark}
In \cite{Hastie22}, the boundedness assumption on $\bSi$ is used primarily to establish equicontinuity in the derivation of the asymptotic limit of the bias term
$B_X(\hat{\bm\beta};\bm\beta)$.
By contrast, this assumption is not required when deriving the asymptotic limit of the variance term
$V_X(\hat{\bm\beta};\bm\beta)$.
In the present paper, we remove the boundedness assumption by decomposing the bias term into two components: one associated with the spiked eigenvalues and the other associated with the bulk eigenvalues.

When the population covariance matrix $\bSi$ has the spiked structure described above, \cite{hu26} shows that the $k$th largest sample eigenvalue associated with the population spike $\al_k$ has the asymptotic limit
\[
\phi(\al_k)
=
\al_k
\left(
1+\frac{\gamma}{\al_k-1}
\right).
\]
\end{remark}
\begin{remark}
	In studying the limiting behavior of the generalization error under a spiked covariance structure, \cite{M19} considered a fixed number of diverging spikes, whereas \cite{li24,li25} focused on a population covariance structure with a single spike. In contrast, our framework allows the number of spikes $M$ either to remain fixed or to diverge, subject to the growth condition
	\[
	\frac{M}{n^{1/4}}\to 0.
	\]
\end{remark}	

\begin{assumption}\label{ass3}
	As $\min\{p,n\}\to\infty$, $H_n:=F^{\bSi}\stackrel{d}{\rightarrow}H$, where $H$ is a distribution function on the real line.
\end{assumption} 
According to \cite{Silverstein95S}, under Assumptions \ref{ass2}--\ref{ass3}, we have $F^{\hat{\bSi}}\stackrel{d}{\rightarrow}F^{\gamma,H}$ almost surely, where $F^{\gamma,H}$ is a nonrandom distribution function whose Stieltjes transform $m:=m_{F^{\gamma,H}}(z)  $ satisfies the following equation:
\begin{align}\label{Sequation}
	m=\int \frac{1}{t(1-\gamma-\gamma z m)-z} d H(t).
\end{align}
In the sequel, we call $F^{\gamma,H}$ the limiting spectral distribution (LSD) of the sample covariance matrix $\hat{\bSi}$.

\section{Main results}	\label{mainresult}
\subsection{Limiting risk}
We are now in a position to state our main theorems. The following results characterize the asymptotic limits of the prediction risk under a generalized spiked covariance model in both the underparameterized and overparameterized regimes.
\begin{thm}\label{thm1}
	Suppose that the observations follow the model specified in (\ref{model1}) and (\ref{model2}), and assume that $\|\bm\beta\|^2=r^2$ for all $n,p$. Under Assumptions \ref{ass2}--\ref{ass3},  for the min-norm least squares estimator $\bm{\hat{\beta}}$,
	\begin{itemize}
	\item when $p/n\rightarrow \gamma<1$, it holds almost surely that
	\begin{align}\label{riskwhengammaxiaoyu1}
		\lim\limits_{n\rightarrow\infty}R_X(\bm{\hat{\beta}};\bm{\beta})=\sigma^2\frac{\gamma}{1-\gamma};
	\end{align}
	\item when $p/n\rightarrow \gamma\in (1,\infty)$, it holds almost surely that
	\begin{align}
		\label{mean}	B_X(\bm{\hat{\beta}};\bm{\beta})&-(1-\frac{1}{\gamma})r^2-
		\sum_{i=1}^{M}(\al_i-1)\left((1-\frac{1}{\gamma})^2\left\langle \bm u_i,\bm{\beta}\right\rangle^2+O_p(\frac{1}{p})+O_p(\frac{\left\langle \bm u_i,\bm{\beta}\right\rangle}{\sqrt{p}})\right)\rightarrow 0
		,\\
		\label{var}	V_X(\bm{\hat{\beta}};\bm{\beta})&-\sigma^2\frac{1}{\gamma-1} (\sum_{j=1}^{M}\al_j\frac{1}{p}+\frac{p-M}{p})\rightarrow0.
	\end{align}
	Therefore, 
	\begin{align}
		R_X(\bm{\hat{\beta}};\bm\beta)-&\left[\sum_{i=1}^{M}(\al_i-1)\left((1-\frac{1}{\gamma})^2\left\langle \bm u_i,\bm{\beta}\right\rangle^2+O_p(\frac{1}{p})+O_p(\frac{\left\langle \bm u_i,\bm{\beta}\right\rangle}{\sqrt{p}})\right) \notag\right.
		\\
		\phantom{=\;\;}&
		\left.+(1-\frac{1}{\gamma})r^2+\sigma^2\frac{1}{\gamma-1} (\sum_{j=1}^{M}\al_j\frac{1}{p}+\frac{p-M}{p})\right]\rightarrow0, \label{risk}
	\end{align}
	as $n,p\rightarrow\infty$.
	\end{itemize}
\end{thm}		
	
\begin{remark}
	Notably, for $\gamma < 1$, the limiting risk calculated in this paper is consistent with that of \cite{Hastie22}. Below, we offer an explanation from a theoretical standpoint.
	When $p<n$, the variance term $$\frac{\sigma^2}{n}\mathrm{tr}(\hat{\bSi}^{+}\bSi)=\frac{\sigma^2}{n}\mathrm{tr}(\hat{\bSi}^{-1}\bSi)=\frac{\sigma^2}{n}\mathrm{tr}(n\bSi^{-\frac{1}{2}}(\bbZ^{\top}\bbZ)^{-1}\bSi^{-\frac{1}{2}}\bSi)=\frac{\sigma^2}{n}\mathrm{tr}(\frac{\bbZ^{\top}\bbZ}{n})^{-1},$$
	therefore it has the same limit as $\bSi=\bbI_p$, then the remaining proof is provided in \cite{Hastie22}.
\end{remark}	
\begin{remark}
	 For $\gamma>1$, when $\bSi=\bbI_p$, equivalently,
	 $\al_1=\al_2=\cdots=\al_M=1$, the formulas in
	 \eqref{mean} and \eqref{var} reduce to the corresponding result in
	 Theorem~1 of \cite{Hastie22}. If
	 \[
	 \sum_{j=1}^{M}\al_j=o(p),
	 \]
	 then the limiting variance satisfies
	 \[
	 V_X(\hat{\bm\beta};\bm\beta)
	 \rightarrow
	 \frac{\sigma^2}{\gamma-1},
	 \]
	 which agrees with the corresponding asymptotic result in
	 \cite{Hastie22}. By contrast, when
	 \[
	 \sum_{j=1}^{M}\al_j=\Omega(p),
	 \]
	 the limiting variance is affected by the spiked eigenvalues. Thus, when the aggregate spike strength is negligible relative to $p$, only the asymptotic bias component is modified by the spiked covariance structure. When the aggregate spike strength is of order $p$ or larger, both the bias and variance components of the asymptotic prediction risk may be affected.	 
\end{remark}

Moreover, it follows from \eqref{risk} that the alignment between the regression coefficient $\bm\beta$ and the population spike eigenvectors $\bm u_i$ also contributes to the asymptotic generalization error. Following \cite{li25}, rather than estimating the unknown spike directions, the regression coefficient, or the angles between them, we treat the degree of target--spike alignment as given and determine whether it is beneficial or detrimental according to the sign of the corresponding coefficient.

Because the coefficient of
$\langle \bm u_i,\bm\beta\rangle^2$
is positive, a larger value of
$\langle \bm u_i,\bm\beta\rangle^2$
leads to a larger asymptotic prediction risk. Thus, stronger alignment between $\bm\beta$ and the corresponding population spike direction is detrimental to generalization. Since $\|\bm\beta\|_2^2=r^2$ and $\|\bm u_i\|_2=1$, we have
\begin{align*}
	\langle \bm u_i,\bm\beta\rangle^2
	&=
	\|\bm\beta\|_2^2
	\|\bm u_i\|_2^2
	\cos^2\!\angle(\bm\beta,\bm u_i) \\
	&=
	r^2\cos^2\!\angle(\bm\beta,\bm u_i)
	\in[0,r^2].
\end{align*}
To some extent, the prediction risk can be controlled by adjusting the degree of alignment.
If one wishes to control the prediction risk and find a more suitable interpolating estimator, we can search in directions that are orthogonal to all spike eigenvectors.

For ease of illustration, we provide a special case when $M=1$ in the following example.
\begin{example}\label{example}
	Under the same assumptions as Theorem \ref{thm1}, and we assume $M=1$, then we have
	\begin{itemize}
	\item when $\gamma<1$,
	\begin{align*}
		\lim\limits_{n\rightarrow\infty}R_X(\bm{\hat{\beta}};\bm{\beta})=\sigma^2\frac{\gamma}{1-\gamma};
	\end{align*}
	\item when $\gamma>1$,
	\begin{align*}
		R_X(\bm{\hat{\beta}};\bm{\beta})-\left[\sigma^2\frac{1}{\gamma-1}(\frac{\al_1}{p}+1)+(1-\frac{1}{\gamma})r^2+(\al_1-1)(1-\frac{1}{\gamma})^2\left\langle\bm u_1,\bm\beta \right\rangle^2 \right]\rightarrow0.
	\end{align*}
	\end{itemize}
\end{example}
Next, we compare the theoretical results for the asymptotic risk with the special case where $\bSi=\bbI_p$. We mainly compare the theoretical results of Theorem \ref{thm1} with those of Theorem 1 in \cite{Hastie22}. For clarity, we summarize the comparison in Table \ref{tab:risk_comparison}.
\begin{table}[htbp]
	\centering
	\caption{Comparison of the asymptotic risk when $\bSi$ is a spiked model and $\bSi=\bbI_p$.}
	\label{tab:risk_comparison}
	\begin{tabular}{p{3cm} p{5cm} p{5cm}}
		\toprule
		& Variance term & Bias term \\
		\midrule
		$M$ is fixed or $M\rightarrow\infty$, $M/n^{1/4}\rightarrow0$ & When $\sum_{i=1}^{M}\al_i=o(p)$, asymptotic limit of the variance term in risk is the same as $\bSi=\bbI_p$; when it increase to $\sum_{i=1}^{M}\al_i=\Omega(p)$, the spiked eigenvalues contribute to the asymptotic limit of the variance component of the risk. & As long as $\bm\beta$ is not orthogonal to the eigenvectors associated with the spiked eigenvalues, the spiked eigenvalues will affect the bias term in the risk. \\
		\midrule
		Main Intuition &  The spiked eigenvalues contribute to the asymptotic limit of the variance component is of order $O(\sum_{i=1}^{M}\al_i/p)$ & The spiked eigenvalues contribute to the asymptotic limit of the bias component is of order $O(\sum_{i=1}^{M}\al_i)$ \\
		\bottomrule
	\end{tabular}
\end{table}

In Corollary \ref{diagonal case}, we compute closed-form expressions in case $\bSi$ is a diagonal spiked covariance matrix.
\begin{cor}\label{diagonal case}
Let $\bSi=\bgL$. Under Assumptions \ref{ass2}--\ref{ass3}, 
\begin{itemize}
\item when $\gamma<1$, the asymptotic limits of the prediction risk is the same as (\ref{riskwhengammaxiaoyu1});
\item  when $\gamma>1$,  the variance term $V_X(\hat{\bm\beta},\bm\beta)$ is the same as (\ref{var}). For bias term $B_X(\bm{\hat{\beta}};\bm{\beta})$, we have 
\begin{align*}
	B_X(\bm{\hat{\beta}};\bm{\beta})-
	\sum_{i=1}^{M}(\al_i-1)((1-\frac{1}{\gamma})^2\left\langle \bm e_i,\bm{\beta}\right\rangle^2+O_p(\frac{1}{p})+O_p(\frac{\left\langle \bm e_i,\bm{\beta}\right\rangle}{\sqrt{p}}))\rightarrow 0
	.
\end{align*}
\end{itemize}
\end{cor}

The above results show that the asymptotic limits of the prediction risk depends not only on the values of the spiked eigenvalues and the parameter $\gamma$, but also on the components of $\bm\beta$ along the first $M$ directions.

Corollary \ref{Mfixedcase} provides the results for a special case when $M$ is fixed.  Here we assume $\bSi$ is a generalized spiked model with structure in Assumption \ref{ass1}.
\begin{cor}\label{Mfixedcase}
Under Assumptions \ref{ass2}--\ref{ass3} and we assume $M$ is fixed,  
\begin{itemize}
\item when $\gamma<1$,
the asymptotic limits of the prediction risk is the same as (\ref{riskwhengammaxiaoyu1}); 
\item when $\gamma>1$, the asymptotic limits of the bias term
 $B_X(\hat{\bm\beta},\bm\beta)$ is the same as (\ref{mean}).
 For the variance term, when 
\begin{itemize}
	\item  $\{\al_j\}_{j=1,\dots,M}$ is bounded, or $\{\al_j\}_{j=1,\dots,M}$ is divergent with $n$ and are of order $o(p)$, we have  
	\begin{align*}
	V_X(\hat{\bm\beta},\bm\beta)-\sigma^2\frac{1}{\gamma-1} \rightarrow0.	
	\end{align*}
	\item $\{\al_j\}_{j=1,\dots,M}$ is divergent with $n$ and are of order $\Omega(p)$, we have 
	\begin{align*}
		V_X(\hat{\bm\beta},\bm\beta)-\sigma^2\frac{1}{\gamma-1}(\frac{\sum_{j=1}^{M}\al_j}{p}+1) \rightarrow0.	
	\end{align*}
\end{itemize}	
\end{itemize}
\end{cor}
Next, we present  the general result of Theorem \ref{thm1} when $M$ is divergent.
\begin{cor}\label{Mdivergentcase}
	Under Assumptions \ref{ass2}--\ref{ass3} and we assume $M$ is divergent and $M=o(n^{1/4})$,  
	\begin{itemize}
	\item when $\gamma<1$,
	the asymptotic limits of the prediction risk is the same as (\ref{riskwhengammaxiaoyu1}); 
	\item when $\gamma>1$,
	the asymptotic limits of $B_X(\hat{\bm\beta},\bm\beta)$  is the same as (\ref{mean}). For variance term, when 
	\begin{itemize}
		\item  $\{\al_j\}_{j=1,\dots,M}$ is bounded, or $\{\al_j\}_{j=1,\dots,M}$ is divergent with $n$ and $\sum_{j=1}^{M}\al_j$ are of order $o(p)$, we have  
		\begin{align*}
			V_X(\hat{\bm\beta},\bm\beta)-\sigma^2\frac{1}{\gamma-1}\frac{p-M}{p} \rightarrow0.	
		\end{align*}
		\item $\{\al_j\}_{j=1,\dots,M}$ is divergent with $n$ and $\sum_{j=1}^{M}\al_j$ are of order $\Omega(p)$, we have $V_X(\hat{\bm\beta},\bm\beta)$ has the same form as (\ref{var}).
\end{itemize}	
\end{itemize}
\end{cor}

\subsection{Risk analysis}\label{simulation}
\begin{figure}[]
	\centering
	\subfigure[Risk curves when $\bSi=\bbI_p$]{
		\includegraphics[width=8.7cm,height=6cm]{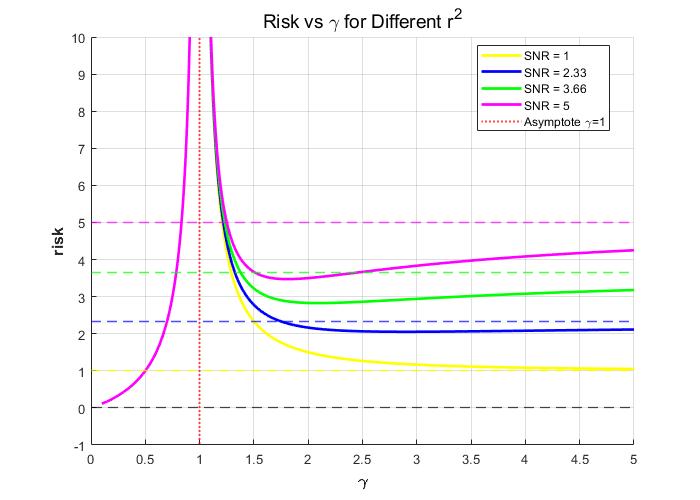}\label{m1}}
	\subfigure[Risk curves when $\bSi$ is a spiked model with $M=5$ and $\sum_{i=1}^{5}\al_i=50$]{
		\includegraphics[width=8.7cm,height=6cm]{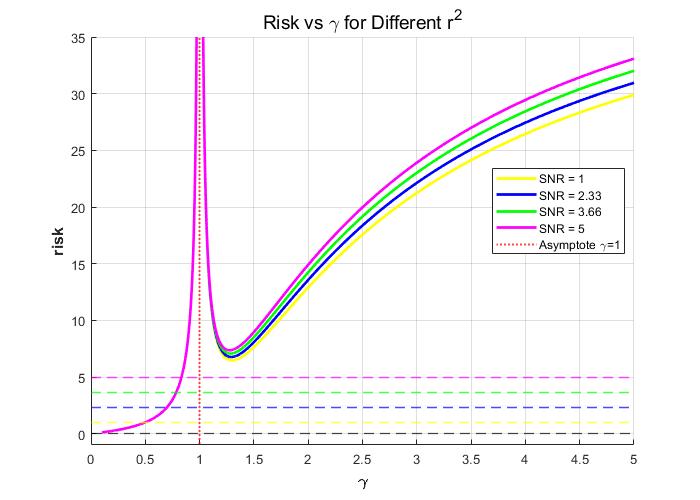}\label{m2}}
	\subfigure[Risk curves when $\bSi$ is a spiked model with $M=5$ and $\sum_{i=1}^{5}\al_i=200$]{
		\includegraphics[width=8.7cm,height=6cm]{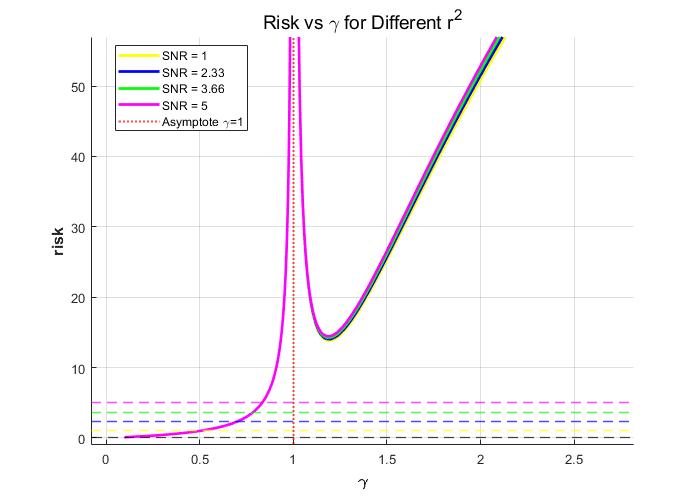}\label{m3}}
	\caption{The asymptotic risk curves as a function of the aspect ratio $\gamma$. Here, we provide a comparison when  $\bSi$ is an identity matrix (Figure \ref{m1}) and the scenario where $\bSi$ is a spiked covariance model with 5 spikes and $\sum_{i=1}^{5}\al_i=50$ (Figure \ref{m2}), and $\sum_{i=1}^{5}\al_i=200$ (Figure \ref{m3}). $\mathrm{SNR}$ varies from 1 to 5, and $\sigma^2=1$. 
	We set $p=300$, $n=p/\gamma$, and $ \left\langle \bm{u_i},\bm   \beta\right\rangle=1  $, where $i=1,2,\dots,5$. The null risk curves corresponding to different $\mathrm{SNRs}$ are denoted by horizontal dashed lines of varying colors in the figure. The specific choices of the $\mathrm{SNRs}$ and $\sigma^2$ in the plot are selected with reference to the study by \cite{Hastie22}. Figure \ref{m1} is also consistent with Figure 2 in \cite{Hastie22}.
	}
	\label{figduibi}
\end{figure}

\begin{figure}[htbp]
	\centering
	\includegraphics[width=16cm,height=8cm]{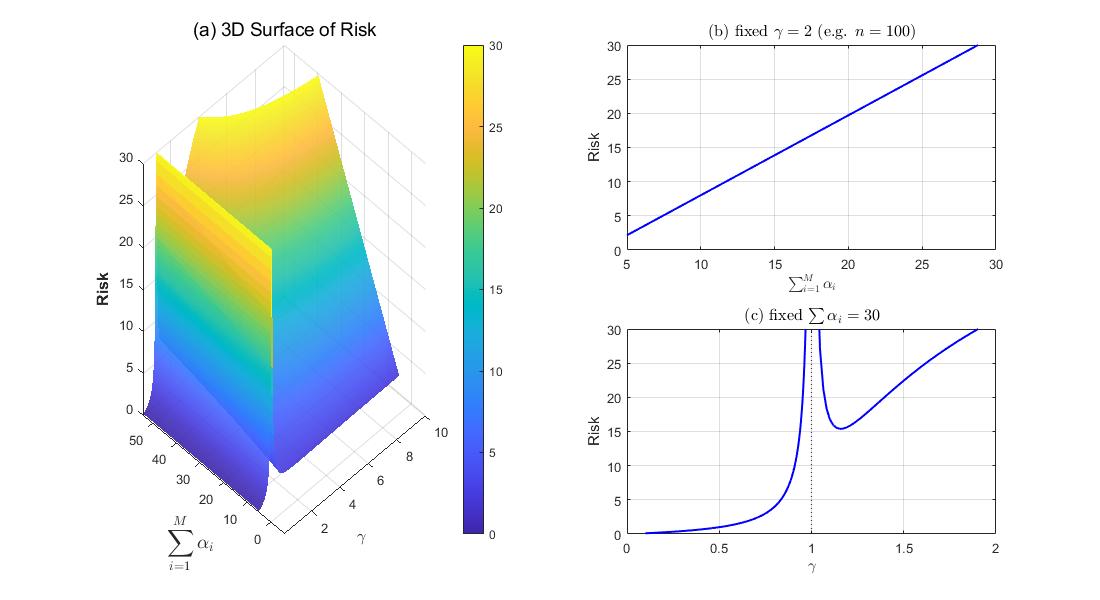}\label{}
	\caption{Left: The asymptotic risk surface varying with $\gamma$ and the summation of spiked eigenvalues $\sum_{i=1}^{M}\al_i$. Right: The corresponding cross-sections of the left surface at $\gamma = 2$ and at  $\sum_{i=1}^{M}\al_i=30$. We set $\mathrm{SNR}=2.33$, $\sigma^2=1$, $p=200$, $M=5$, and $ \left\langle \bm{u_i},\bm   \beta\right\rangle=1  $, where $i=1,2,\dots,5$.
	}
	\label{fig_risk_3D}
\end{figure}

\begin{figure}[htbp]
	\centering
	\includegraphics[width=15cm,height=7cm]{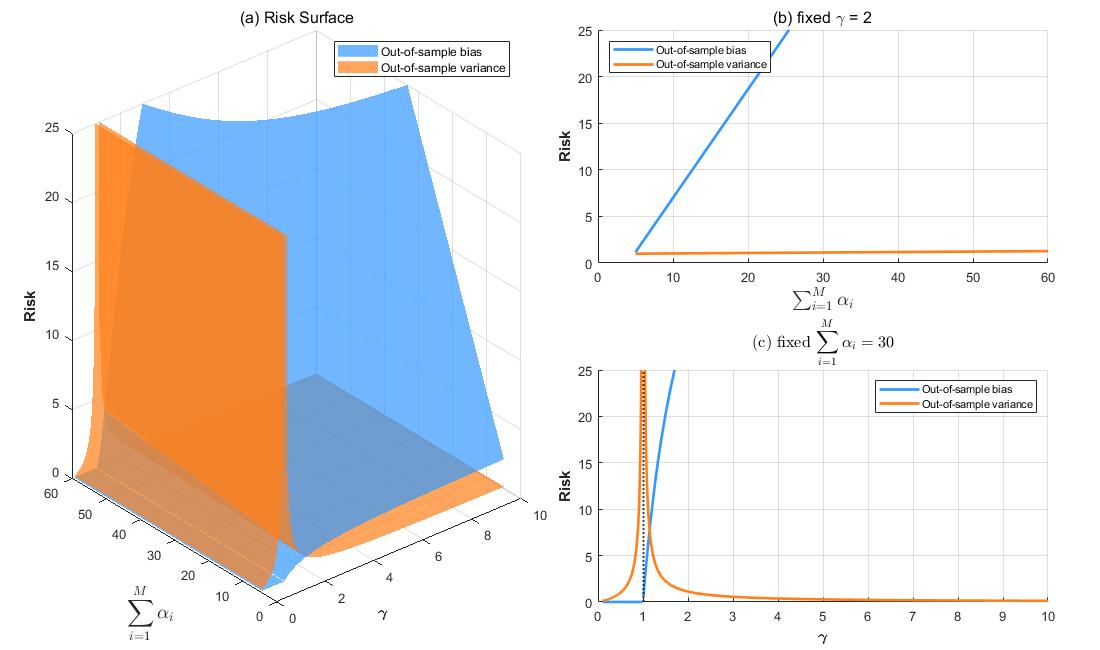}
	\caption{Left: Asymptotic risk surfaces for the bias (blue) and variance (orange) terms as functions of $\gamma$ and the summation of spiked eigenvalues $\sum_{i=1}^{M}\al_i$. Right: Cross-sectional profiles of the left surfaces at $\gamma = 2$ and at  $\sum_{i=1}^{M}\al_i=30$. We set $\mathrm{SNR}=2.33$, $\sigma^2=1$, $p=200$, $M=5$ and $ \left\langle \bm{u_i},\bm   \beta\right\rangle=1  $, where $i=1,2,\dots,5$.
	}
	\label{fig_risk_decomposition_3D}
\end{figure}

\begin{figure}[htbp]
	\centering
	\includegraphics[width=15cm,height=7cm]{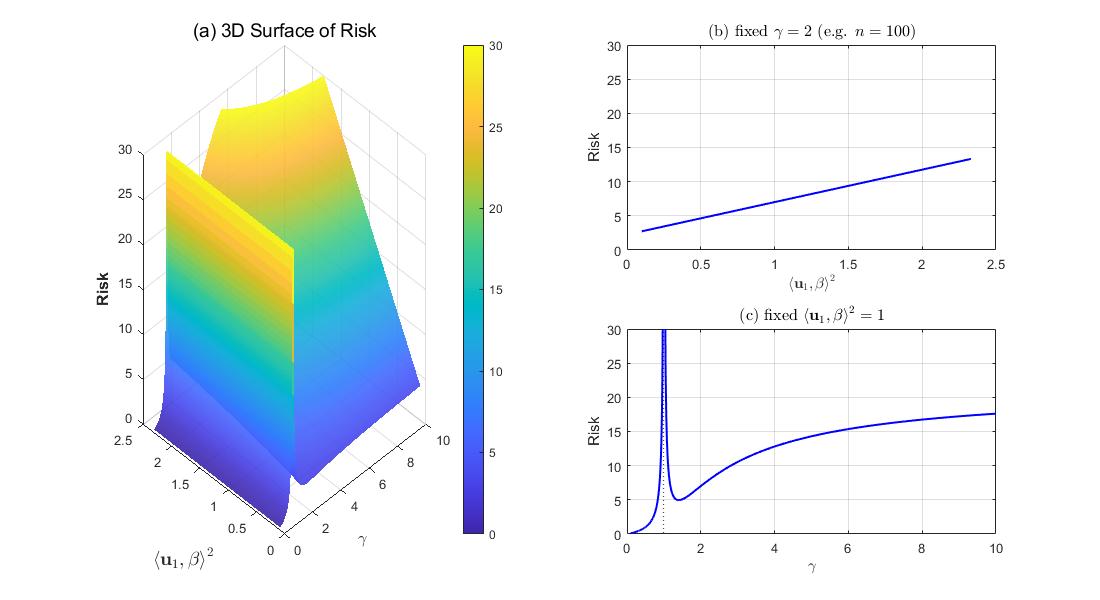}
	\caption{Left: The asymptotic risk surface varying with $\gamma$ and $\left\langle \bm{u_i},\bm   \beta\right\rangle^2$. Right: The corresponding cross-sections of the left surface at $\gamma = 2$ and at  $\left\langle \bm{u_i},\bm   \beta\right\rangle^2=1$. We set $\mathrm{SNR}=2.33$, $\sigma^2=1$, $p=200$, $M=1$, and $\al_1=20$.
	}
	\label{fig_risk_innerproduct_3D}
\end{figure}

In this section, we provide an analysis of the results presented in Theorem \ref{thm1} and conducted some simulations (Figures \ref{figduibi}--\ref{fig_risk_innerproduct_3D}). We denote the signal-to-noise ratio as $\mathrm{SNR}=\|  \bm\beta\|_2^2/\sigma^2$.  Let $\mathrm{SNR}=r^2 / \sigma^2$. 
From analysis before Example \ref{example}, we know that the range of
$\langle \bm u_i,\bm\beta\rangle^2$
is $[0,r^2]$. Therefore, in Figures \ref{figduibi}--\ref{fig_risk_innerproduct_3D} we set $\langle \bm u_i,\bm\beta\rangle^2=1$.
For the sake of clarity, we first give a description of Figures \ref{figduibi}--\ref{fig_risk_innerproduct_3D} as follows:
\begin{itemize}
	\item Figure \ref{figduibi}  illustrates the risk as a function of $\gamma$ under three different population structures.
	\item Figure \ref{fig_risk_3D} presents a surface plot of the risk as a function of the variables $\gamma$ and the summation of spiked eigenvalues $\sum_{i=1}^{M}\al_i$.
	
	\item Figure \ref{fig_risk_decomposition_3D}  presents surface plots of the bias and variance components of the risk as functions of the variables $\gamma$ and  $\sum_{i=1}^{M}\al_i$.
	\item Figure \ref{fig_risk_innerproduct_3D} presents a surface plot of the risk as a function of the variables $\gamma$ and $\left\langle \bm{u_i},\bm   \beta\right\rangle^2$. 
\end{itemize}

The following facts are immediate from the risk in (\ref{riskwhengammaxiaoyu1}) and (\ref{risk}). See Figure \ref{figduibi} for an accompanying plot when $\mathrm{SNR}$ varies from 1 to 5. 
\begin{itemize}
	\item For $\gamma\in(0,1)$,  there is no bias and then 
	the risk equals the variance. Since risk is only related to $\gamma$ and $\sigma^2$ when $\gamma\in(0,1)$, consequently, the curves in the Figures \ref{m1} and \ref{m2} coincide regardless of the specific $\mathrm{SNR}$ values. The risk is increases with $\gamma$. For $\gamma >1$, the risk is constituted by the bias term and the variance term. From (\ref{risk}) we know that the bias term increases with $\gamma$ and the variance term decreases with $\gamma$.
	\item By comparing Figures \ref{m1} and \ref{m2}, it is apparent that when $\gamma>1$, the presence of the spikes accelerates the convergence of the risk curve to its minimum, followed by a rapid rebound. This suggests that the bias term plays a more dominant role in this scenario. It is consistent with analyses in the Table \ref{tab:risk_comparison}, since the spiked eigenvalues contribute to the asymptotic limit of the variance component is of order $O(\sum_{i=1}^{M}\al_i/p)$, while to the bias term is of order $O(\sum_{i=1}^{M}\al_i)$.
	Moreover, on $\gamma\in (1,\infty)$, the null risk beats risk curves when $\mathrm{SNR}$ varies from 1 to 5. 
	A comparison of the curves at the same SNR across Figures \ref{m1} and \ref{m2} reveals that the presence of spiked eigenvalues causes the risk curve to attain its minimum earlier and, relative to the isotropic case, leads to a higher minimum risk.
	It is also observed that as the $\mathrm{SNR}$ decreases, the position of the minimum increases, whereas the corresponding  minimum risk decreases. Both $\sum_{i=1}^{M}\al_i$ and $\mathrm{SNR}$ affects risk. From Figure \ref{m2} we can conclude that a higher $\mathrm{SNR}$ accelerates the divergence of the risk curve.
	\item By comparing Figures \ref{m2} and \ref{m3},  when $\gamma>1$, regardless of the SNR considered in the figure, a consistent conclusion is that larger spiked eigenvalues lead to a higher minimum of the risk curve.
	
\end{itemize}	

We provide the asymptotic risk surface varying with $\gamma$ and the summation of spiked eigenvalues $\sum_{i=1}^{M}\al_i$ in Figure \ref{fig_risk_3D}. Compared with Figure \ref{figduibi}, it introduces an additional dimension, $\sum_{i=1}^{M}\alpha_i$, which allows us to observe the trend of the risk with respect to $\sum_{i=1}^{M}\alpha_i$ as well. When $\gamma < 1$ and $\sum_{i=1}^{M}\alpha_i$ is bounded, the asymptotic risk is not affected by the value of $\sum_{i=1}^{M}\alpha_i$, which aligns with our theoretical findings (Theorem \ref{thm1}). When $\gamma > 1$, the risk universally exhibits a U-shaped trend, first decreasing and then increasing. As $\sum_{i=1}^{M}\alpha_i$ grows, the minimum of the risk curve emerges earlier, while its minimum value becomes progressively larger. For a fixed $\gamma$ and $\gamma > 1$, the asymptotic risk grows linearly with the increase of $\sum_{i=1}^{M}\alpha_i$. Figure \ref{fig_risk_decomposition_3D} provides the bias-variance decomposition of the asymptotic risk presented in Figure \ref{fig_risk_3D}. It can be observed that when $\gamma > 1$, the bias term plays a dominant role in determining the behavior of the asymptotic risk.

Figure \ref{fig_risk_innerproduct_3D} shows the asymptotic risk surface varying with $\gamma$ and $\left\langle \bm{u_i},\bm   \beta\right\rangle^2$. For ease of illustration, we focus on the case $(M=1)$ to examine the effect of $\left\langle \bm{u_i},\bm   \beta\right\rangle^2$ on the risk. Here $\left\langle \bm{u_i},\bm   \beta\right\rangle^2$ denotes the alignment of $\bm\beta$ and the associated eigenvector $ \bm u_1 $ of spiked eigenvalues $\al_1$. Since 
\begin{align*}
	\left\langle \bm{u_i},\bm   \beta\right\rangle^2\in [0,r^2],
\end{align*}
and we set $\mathrm{SNR}=r^2/\sigma^2=2.33$, therefore $r^2=2.33$. Accordingly, in Figure \ref{fig_risk_innerproduct_3D}, the y-axis ranges from $0$ to $2.33$. A value of $0$ indicates that $\bm\beta$ is orthogonal to $\bm u_1$, while a larger value on the y-axis corresponds to a smaller angle between them. From Figure \ref{fig_risk_innerproduct_3D}, it is observed that the risk depends linearly on $\left\langle \bm{u_i},\bm   \beta\right\rangle^2$. Moreover, the less aligned $\bm u_1$ and $\bm\beta$ are, the lower the risk becomes, indicating that, within our framework, their alignment is detrimental to prediction performance.

\section{Statistical implications}\label{stat imp}
\subsection{Benign, tempered, and catastrophic overfitting under spiked covariance model}
In overparameterized regression, interpolation does not necessarily imply poor generalization. The phenomenon of benign overfitting shows that a minimum-norm interpolator can achieve vanishing excess risk despite perfectly fitting noisy training data (see \cite{Bartlett20}). However, existing characterizations have largely focused on isotropic or regular covariance structures. In the presence of spikes, the interaction between spike strength, the number of spikes, and signal alignment can fundamentally alter the generalization behavior. Unlike \cite{li25}, which focuses on a single--spike covariance model, we consider a more general setting with multiple spikes. In this framework, the number of spikes $M$ becomes an additional structural parameter that plays a crucial role in determining the generalization behavior of interpolating estimators.
In the following, we provide the definitions of the excess risk and the taxonomy of overfitting.
\begin{proposition}\label{prop}
	Under the spiked model defined in Assumption \ref{ass1}, the excess risk is defined as
	\begin{align*}
	R_\gamma\triangleq	\lim\limits_{n\rightarrow\infty,p/n\rightarrow\gamma}R_X(\bm{\hat{\beta}};\bm\beta)-r^2,
	\end{align*}
	where $r^2$ is the risk of the null estimator, represents the irreducible noise.
	
	Then following \cite{Bartlett20} and \cite{Mallinar22}, the minimum--norm interpolator exhibits three possible asymptotic behaviors:
	\begin{itemize}
		\item Benign overfitting:
		\begin{align*}
		\text{if}~~	\lim\limits_{\gamma\rightarrow\infty}R_\gamma=0;
		\end{align*}
		\item Tempered overfitting:
		\begin{align*}
		\text{if}~~	0<\lim\limits_{\gamma\rightarrow\infty}R_\gamma<\infty;
		\end{align*}
		\item Catastrophic overfitting:
		\begin{align*}
		\text{if}~~	\lim\limits_{\gamma\rightarrow\infty}R_\gamma=\infty.
		\end{align*}
	\end{itemize}
\end{proposition}
The following theorem discusses how the spike strength, the number of spikes, overparameterization $(\gamma = p/n)$ and signal alignment affect the overfitting phase classification.
\begin{thm}\label{thm_overfitting}
		Assume the model (\ref{model1}) and (\ref{model2}). Also assume that $\|\bm\beta\|^2=r^2$ for all $n,p$. Under Assumptions \ref{ass2}--\ref{ass3}, depending on whether $\bm\beta$ is orthogonal to $\{ \bm u_j \}_{j=1}^{M}$,  we consider the following case:
		\begin{itemize}
			\item when $\bm\beta$ is orthogonal with $\{ \bm u_j \}_{j=1}^{M}$, and 
			
			(i) $\sum_{j=1}^{M}\al_j=o(p)$, or $\sum_{j=1}^{M}\al_j=\Omega(p)$ but the divergence rate of $\sum_{j=1}^{M}\al_j/p$ is slower than the rate at which $\gamma\rightarrow\infty$, then benign overfitting arises;
			
			(ii) the divergence rate of $\sum_{j=1}^{M}\al_j/p$ is of the same order as the rate at which $\gamma\rightarrow\infty$, then tempered overfitting arises;
			
			(iii) the divergence rate of $\sum_{j=1}^{M}\al_j/p$ is faster than the rate at which $\gamma\rightarrow\infty$, then catastrophic overfitting arises;
			\item when $\bm\beta$ is not orthogonal to all vectors in set $\{ \bm u_j \}_{j=1}^{M}$, or it is not orthogonal to any vector in set $\{ \bm u_j \}_{j=1}^{M}$. Let $\{ \bm u_j \}_{j\in J_k}$ denote the part that is not orthogonal to $\bm\beta$, and
			
			(i) $\sum_{j\in J_k}\al_j$ is of constant order, and the divergence rate of $\sum_{j=1}^{M}\al_j/p$ is not faster than the rate at which $\gamma\rightarrow\infty$, then tempered overfitting arises;
			
			(ii) $\sum_{j\in J_k}\al_j$ is of constant order, and the divergence rate of $\sum_{j=1}^{M}\al_j/p$ is faster than the rate at which $\gamma\rightarrow\infty$, or $\sum_{j\in J_k}\al_j$ is divergent, then catastrophic overfitting arises.
		\end{itemize}
\end{thm}

Similar to \cite{li25}, we summarize the results of Theorem \ref{thm_overfitting} in Table \ref{fig_overfitting}.

\begin{table}[htbp]
	\centering
	\caption{Asymptotic Generalization Regimes. This table summarizes conditions for when overfitting is benign, tempered, or catastrophic in the limit where $\gamma=p/n$ and subsequently $\gamma\rightarrow\infty$. The behavior depends on the spike strength, the number of spikes, overparameterization $\gamma$ and signal alignment}
	\label{fig_overfitting}
	\begin{tabular}{p{3cm} p{3cm} p{3cm} p{3cm}}
		\toprule
		Regime &Benign & Tempered  & Catastrophic \\
		\midrule
		Case 1($\bm\beta$ is orthogonal with $\{ \bm u_j \}_{j=1}^{M}$): &  $\sum_{j=1}^{M}\al_j=o(p)$, or $\sum_{j=1}^{M}\al_j=\Omega(p)$ but the divergence rate of $\sum_{j=1}^{M}\al_j/p$ is slower than the rate at which $\gamma\rightarrow\infty$ &  the divergence rate of $\sum_{j=1}^{M}\al_j/p$ is of the same order as the rate at which $\gamma\rightarrow\infty$ & the divergence rate of $\sum_{j=1}^{M}\al_j/p$ is faster than the rate at which $\gamma\rightarrow\infty$\\
		\midrule
	Case 2(when $\bm\beta$ is not orthogonal to all vectors in set $\{ \bm u_j \}_{j=1}^{M}$, or it is not orthogonal to any vector in set $\{ \bm u_j \}_{j=1}^{M}$):	& \multicolumn{1}{c}{$\times$} & $\sum_{j\in J_k}\al_j$ is of constant order, and the divergence rate of $\sum_{j=1}^{M}\al_j/p$ is not faster than the rate at which $\gamma\rightarrow\infty$ & $\sum_{j\in J_k}\al_j$ is of constant order, and the divergence rate of $\sum_{j=1}^{M}\al_j/p$ is faster than the rate at which $\gamma\rightarrow\infty$, or $\sum_{j\in J_k}\al_j$ is divergent\\
		\bottomrule
	\end{tabular}
\end{table}

Our results show that benign overfitting is not solely a consequence of overparameterization. Instead, it emerges from a delicate balance between spike strength, the number of spikes, overparameterization $(\gamma = p/n)$ and target alignment. Spiked covariance structures therefore provide a principled framework for understanding when interpolation is statistically harmless and when it becomes unstable.

\subsection{Double descent under spiked covariance model}
Since the effect of the spike eigenvalues on the prediction risk is manifested only when $\gamma > 1$, we restrict our attention to the $\gamma > 1$ case below. Based on our earlier analysis, the effect of spikes on the risk is not monotonic. Instead, it is determined by the strength and quantity of the spikes, along with the alignment between the regression vector $\bm\beta$ and the spike eigenvectors $\{\bm u_i\}_{i=1}^M$. From Figure \ref{figduibi} we can conclude that when $\left\langle \bm\beta,\bm u_i\right\rangle\neq0 $, the spike structure can destroy the benign overparameterization phenomenon of the isotropic model, and it can induce a new failure mode of double descent in the overparameterized regime. 
A natural question arises: does the same conclusion hold when $\left\langle \bm\beta,\bm u_i\right\rangle=0 $? We answer this question below.

When $\left\langle \bm\beta,\bm u_i\right\rangle=0 $, then the asymptotic risk in Theorem \ref{thm1} is 
\begin{align*}
	R_X(\bm{\hat{\beta}};\bm\beta)-\left[(1-\frac{1}{\gamma})r^2+\sigma^2\frac{1}{\gamma-1}(\sum_{j=1}^{M}\al_j\frac{1}{p}+\frac{p-M}{p})\right]\rightarrow0,
\end{align*}
therefore, when $\sum_{j=1}^{M}\al_j=o(p)$, the asymptotic risk is the same as isotropic case ($\bSi=\bbI_p$); when  $\sum_{j=1}^{M}\al_j=\Omega(p)$, it is only at this point that the impact of the spiked eigenvalues on the prediction risk becomes apparent. As $\gamma$ increases, the variance term decreases while the bias term increases. However, the presence of the spikes causes the variance term to dominate; consequently, the overall risk exhibits a downward trend. This is in stark contrast to our earlier conclusion for the $\left\langle \bm\beta,\bm u_i\right\rangle\neq0$ case. When $\left\langle \bm\beta,\bm u_i\right\rangle\neq0$, the bias term dominates, and thus the risk ultimately exhibits an increasing trend as $\gamma$ grows. This finding offers valuable insights: for real-world data where spikes may emerge (such as financial data), if one wishes to control the prediction risk and find a more suitable interpolating estimator, we can search in directions that are orthogonal to all spike eigenvectors.

\begin{figure}[]
	\centering
	\subfigure[Risk curves when $\bSi$ is a spiked model with $M=5$ and $\sum_{i=1}^{5}\al_i=50$]{
		\includegraphics[width=7cm,height=6cm]{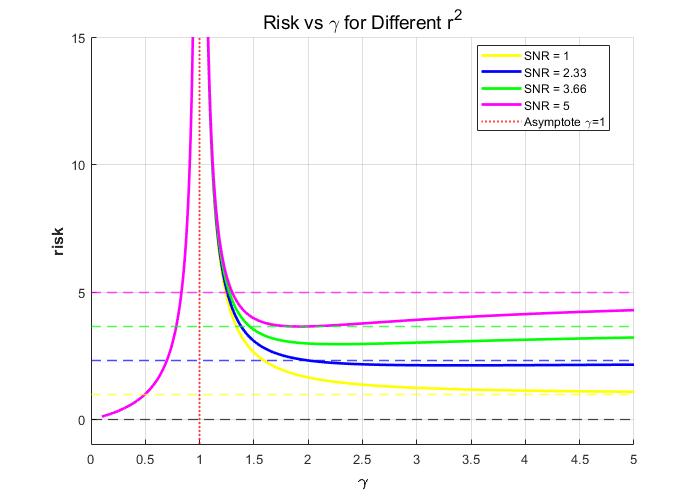}\label{bounded_oro}}
	\subfigure[Risk curves when $\bSi$ is a spiked model with $M=5$ and $\sum_{i=1}^{5}\al_i=500$]{
		\includegraphics[width=7cm,height=6cm]{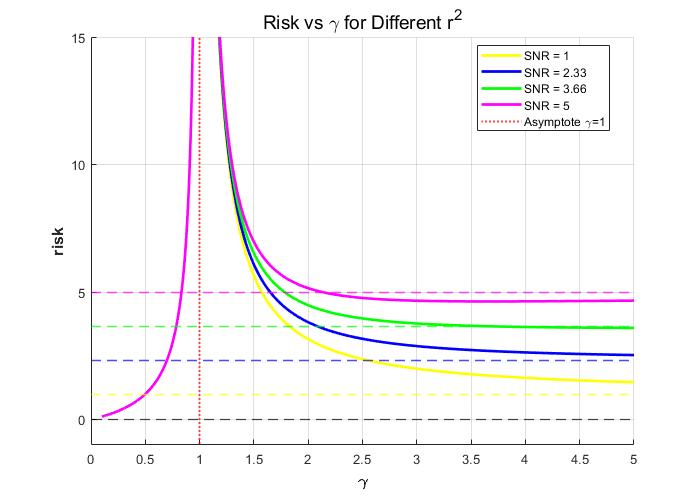}\label{divergent_oro}}
	\caption{The asymptotic risk curves as a function of the aspect ratio $\gamma$. Here, we provide a comparison when  $\bSi$ is a spiked covariance model with 5 spikes and $\sum_{i=1}^{5}\al_i=50$ (Figure \ref{bounded_oro}), and $\sum_{i=1}^{5}\al_i=200$ (Figure \ref{divergent_oro}). $\mathrm{SNR}$ varies from 1 to 5, and $\sigma^2=1$. We set $p=300$, $n=p/\gamma$. Different from Figure \ref{figduibi}, we let $ \left\langle \bm{u_i},\bm   \beta\right\rangle=0  $, where $i=1,2,\dots,5$. The null risk curves corresponding to different $\mathrm{SNRs}$ are denoted by horizontal dashed lines of varying colors in the figure. The specific choices of the $\mathrm{SNRs}$ and $\sigma^2$ in the plot are selected with reference to the study by \cite{Hastie22}. 
	}
	\label{figduibi_oro}
\end{figure}

In the following, to support the preceding theoretical explanation, we present the risk curves for $\left\langle \bm\beta,\bm u_i\right\rangle=0$ in Figure \ref{figduibi_oro}. From Figure \ref{bounded_oro}, we find that when $\sum_{i=1}^{5}=50$ and it can be regarded as being of order $o(p)$. Then  the asymptotic risk is the same as isotropic case ($\bSi=\bbI_p$); from Figure \ref{divergent_oro}, when  $\sum_{j=1}^{M}\al_j=500$, and it can be regarded as being of order $\Omega(p)$,  then overall risk exhibits a downward trend. This is consistent with our earlier analysis.

\section*{Appendix}	\label{appendix}
In this section, we provide some useful lemmas.
\begin{lemma}[\cite{Johnstone18}] 
Under Assumptions \ref{ass2}--\ref{ass3}, the squared inner product between $\nu$-th sample and population eigenvector has the following limit:
\begin{align}
	&\left\langle \bm u_{\nu,} \bm v_\nu\right\rangle^2 \xrightarrow{a.s.}\left\{\begin{array}{lll}
		\frac{1-\gamma /\left(\al_\nu-1\right)^2}{1+\gamma /\left(\al_\nu-1\right)} & \text { if } & \al_\nu>1+\sqrt{\gamma} \\
		0 & \text { if } & \al_\nu \leq 1+\sqrt{\gamma}
	\end{array}\right.
\end{align}
\end{lemma}
\begin{lemma}[\cite{bai12}]
	Assume Assumptions \ref{ass2}--\ref{ass3} hold. Define $\phi(\al_k)=\al_k(1+\frac{\gamma}{\al_k-1})$. Let $\al_k$ be a generalized spiked eigenvalue satisfying  and $\phi^{\prime}(\al_k)>0$ (distant spike), then the sample spiked eigenvalues $\lambda_k$ converge almost surely to $\phi(\al_k)$.
\end{lemma}	
\begin{lemma}[\cite{hu26}]\label{hu26}
 Suppose Assumptions \ref{ass2}--\ref{ass3} in Section \ref{notation} and one of the following two conditions hold:
	\begin{itemize}
		\item[(i)] $z_{i j}$ is Gaussian and $M=o\left(n^{1 / 3}\right)$,
		\item[(ii)] $\mathbb{E} z_{i j}^4<\infty$ and $M=o\left(n^{1 / 4}\right)$.
	\end{itemize}
	Then, for all $1 \leq i \leq M$, we have
	\begin{align*}
		\frac{s_i}{\phi\left(\al_i\right)} \xrightarrow{p} 1, \\
		\frac{-\hat{\underline{S}}_{n,i}^{-1}}{s_i} \xrightarrow{p} 1,
	\end{align*}
	where $\underline{\hat{S}}_{n,i}^{-1}=1 / \hat{\underline{S}}_{n,i}$, and $\hat{\underline{S}}_{n,i}=-\frac{1-\gamma}{s_i}+\frac{1}{n}\sum_{j=M+1}^{p}\frac{1}{s_j-s_i}$.
\end{lemma}
\begin{remark}
	Lemma \ref{hu26} is borrowed in \cite{hu26}. This lemma provides a method for computing the limits to which the sample spiked eigenvalues converge in probability, as well as a method for estimating the corresponding population spiked eigenvalues when $M$ is divergent.
\end{remark}


\begin{thebibliography}{53}
	\providecommand{\natexlab}[1]{#1}
	\providecommand{\url}[1]{\texttt{#1}}
	\expandafter\ifx\csname urlstyle\endcsname\relax
	\providecommand{\doi}[1]{doi: #1}\else
	\providecommand{\doi}{doi: \begingroup \urlstyle{rm}\Url}\fi
	
	\bibitem[Bai and Ng (2002)Bai, Ng]{bai02}
	Jushan Bai and Serena Ng.
	\newblock Determining the number of factors in approximate factor models.
	\newblock \emph{Econometrica}, 70\penalty0 (1):\penalty0
	191--221, 2002.
	
	
	
	
	
	
	\bibitem[Bai and Zhou(2008)]{Bai08}
	Zhidong Bai and Wang Zhou.
	\newblock {Large sample covariance matrices without
		independence structures in columns}.
	\newblock \emph{Statistica Sinica},18\penalty0:\penalty0 425 --
	442, 2008.
	
	
	
	\bibitem[Bai and Yao(2012)]{bai12}
	Zhidong Bai and Jianfeng Yao.
	\newblock On sample eigenvalues in a generalized spiked population model.
	\newblock \emph{Journal of Multivariate Analysis}, 106:\penalty0 167--177,
	2012.
	
	\bibitem[Bai et al.(2018)]{bai18}
	Zhidong Bai,  Kwok Pui Choi and  Yasunori Fujikoshi.
	\newblock Consistency of AIC and BIC in estimating the number of significant components in high-dimensional principal component analysis.
	\newblock \emph{The Annals of Statistics}, 46:\penalty0(3) 1050--1076,
	2018.
	
	\bibitem[Baik and Silverstein(2006)]{baik06}
	Jinho Baik and Jack~W. Silverstein.
	\newblock Eigenvalues of large sample covariance matrices of spiked population
	models.
	\newblock \emph{Journal of Multivariate Analysis}, 97\penalty0 (6):\penalty0
	1382--1408, 2006.
	
	\bibitem[Bartlett et al.(2020)]{Bartlett20}
	Peter L Bartlett, Philip M Long, Gábor Lugosi and Alexander Tsigler. 
	\newblock Benign overfitting in linear regression.
	\newblock \emph{Proceedings of the National Academy of Sciences},  117\penalty0 (48):\penalty0 30063–30070, 2020.
	
	
	\bibitem[Belkin et al.(2019)]{Belkin19}
	Mikhail Belkin, Daniel Hsu, Siyuan Ma and Soumik Mandal.
	\newblock Reconciling modern machine-learning practice and the classical bias–variance trade-off.
	\newblock \emph{Proceedings of the National Academy of Sciences}, 116:\penalty0(32) 15849--15854,
	2019.
	
	\bibitem[Belkin et al.(2020)]{Belkin20}
	Mikhail Belkin, Daniel Hsu and Ji Xu.
	\newblock Two models of double descent for weak features.
	\newblock \emph{SIAM Journal on Mathematics of Data Science}, 2:\penalty0(4) 1167--1180,
	2020.
	
	
	
	
	
	
	

	
	
	
	\bibitem[Jiang and Bai(2021)]{JiangB21G}
	Dandan Jiang and Zhidong Bai.
	\newblock Generalized four moment theorem and an application to clt for spiked
	eigenvalues of high-dimensional covariance matrices.
	\newblock \emph{Bernoulli}, 27\penalty0 (1):\penalty0 274--294,
	2021.
	
	\bibitem[Johnstone(2001)]{Johnstone01}
	Iain~M. Johnstone.
	\newblock {On the distribution of the largest eigenvalue in principal components
		analysis}.
	\newblock \emph{The Annals of Statistics}, 29\penalty0 (2):\penalty0 295 --
	327, 2001.
	
	\bibitem[Johnstone and Nadler(2017)]{Johnstone17}
	Iain~M. Johnstone and Boaz Nadler.
	\newblock Roy's largest root test under rank-one alternatives.
	\newblock \emph{Biometrika}, 104\penalty0 (1):\penalty0 181--193, 2017.
	
	\bibitem[Johnstone and Yang(2018)]{Johnstone18}
	Iain~M. Johnstone and Jeha Yang.
	\newblock Notes on asymptotics of sample eigenstructure for
	spiked covariance models with non-Gaussian data.
	\newblock \emph{arXiv:1810.10427},
	2018.
	
	\bibitem[Hastie et~al.(2022)]{Hastie22}
	Trevor Hastie,  Andrea Montanari, Saharon  Rosset and Ryan J Tibshirani.
	\newblock {Surprises in high-dimensional ridgeless least squares interpolation}.
	\newblock \emph{The Annals of Statistics}, 50\penalty0 (2):\penalty0 949--986, 2022.
	
	\bibitem[Hu et~al.(2026)]{hu26}
	Jianwei Hu, Jingfei Zhang, Jianhua Guo and Ji Zhu.
	\newblock {Limiting laws and consistent estimation criteria for fixed and diverging number of spiked eigenvalues}.
	\newblock \emph{Journal of the American Statistical Association}, 2026.
	
	\bibitem[Li and Sonthalia(2024)]{li24}
	Jiping Li and Rishi Sonthalia.
	\newblock Generalization for least squares regression with simple spiked covariances.
	\newblock \emph{arXiv:2410.13991v1},
	2024.
	
	\bibitem[Li and Sonthalia(2025)]{li25}
	Jiping Li and Rishi Sonthalia.
	\newblock Risk phase transitions in spiked regression: alignment driven
	benign and catastrophic overfitting.
	\newblock \emph{arXiv:2510.01414v1},
	2025.
	
	
	\bibitem[Mahdaviyeh and Naulet(2019)]{M19}
	Yasaman Mahdaviyeh and Zacharie Naulet.
	\newblock Risk of the least squares minimum norm estimator
	under the spike covariance model.
	\newblock \emph{arXiv:1912.13421},
	2019.
	
	\bibitem[Mallinar et al.(2022)]{Mallinar22}
Neil Mallinar, James B. Simon and Amirhesam Abedsoltan.
	\newblock Benign, tempered, or catastrophic:
	a taxonomy of overfitting.
	\newblock \emph{36th Conference on Neural Information Processing Systems (NeurIPS 2022).}
	
	
	
	\bibitem[Paul(2007)]{Paul07}
	Debashis Paul.
	\newblock Asymptotics of sample eigenstructure for a large dimensional spiked
	covariance model.
	\newblock \emph{Statistica Sinica}, 17\penalty0 (4):\penalty0 1617--1642, 2007.
	
	\bibitem[Silverstein(1995)]{Silverstein95S}
	Jack~W. Silverstein.
	\newblock Strong convergence of the empirical distribution of eigenvalues of
	large dimensional random matrices.
	\newblock \emph{Journal of Multivariate Analysis}, 55\penalty0 (2):\penalty0
	331--339, 1995.
	
	
\end{thebibliography}
\end{document}